# Pseudo-Newton method for nonlinear equations

W. Chen

**Present mail address** (as a JSPS Postdoctoral Research Fellow): Apt.4, West 1st floor, Himawari-so, 316-2, Wakasato-kitaichi, Nagano-city, Nagano-ken, 380-0926, JAPAN

Permanent affiliation and mail address: Dr. Wen CHEN, P. O. Box 2-19-201, Jiangshu University of Science & Technology, Zhenjiang City, Jiangsu Province 212013, P. R. China

Present e-mail: chenw@homer.shinshu-u.ac.jp
Permanent email: chenwwhy@hotmail.com

**Abstract**: In order to avoid the evaluation of the Jacobian matrix and its inverse, the present author [2] introduced the pseudo-Jacobian matrix with a general applicability of any nonlinear systems of equations. By using this concept, this paper proposes the pseudo-Newton method.

## 1. Introduction.

The Newton method and its variants are of central importance to compute a variety of nonlinear algebraic equation [1]. The time-consuming effort of this method is the iterative evaluation of function values, Jacobian matrix and its inverse. The size of the practical physical and engineering problems is usually very large. Therefore, the solution procedure of these equations is often very costly. An inordinate amount of computing time and storage even prohibits such calculations.

Considerable research effort has been devoted to the development of some efficient nonlinear algorithms to reduce the cost in the evaluation of the Jacobian matrix and its inverse [1]. Very recently, the present author [2] proposed a new concept of the pseudo-Jacobian matrix for stability analysis of nonlinear initial value problems. The objective of this note is to apply this theorem to derive a simple Newton iterative formula which greatly reduces the computing effort in the evaluation of the Jacobian matrix and its inversion..

## 2. Pseudo-Newton method

The pseudo-Jacobian matrix was applied in this section to form a novel variant of Newton method in the solution of the nonlinear problems. Consider the nonlinear equations

$$f(U) = LU + N^{(2)}(U) + N^{(3)}(U) + b = 0, \qquad (1)$$

where $LU$, $N^{(2)}(U)$ and $N^{(3)}(U)$ represent the linear, quadratic and cubic terms of the system of equations. The above equation can be restated by adding one quadratic and two cubic nonlinear terms in two sides

$$LU + 2N^{(2)}(U) + 3N^{(3)}(U) = b + N^{(2)}(U) + 2N^{(3)}(U) \qquad (2)$$

The pseudo-Jacobian matrix is here defined as

$$2N^{(2)}(U) + 3N^{(3)}(U) = \left\{ \left[ 2N^{(2)}(U) + 3N^{(3)}(U) \right] \left[ \frac{1}{n} \left( U^{\circ(-1)} \right)^T \right] \right\} U$$
$$= \left( wv^T \right) U \qquad (3)$$

$U^{\circ(-1)}$ is the Hadamard inverse of vector U. Let

$$\bar{J}(U) = L + wv^T \qquad (4)$$

It is obvious that $\bar{J}(U)$ is an alternative of original Jacobian matrix $J(U)$ of equation (1). Substituting equations (3) and (4) into equation (2) yields

$$\bar{J}(U)U = b + N^{(2)}(U) + 2N^{(3)}(U). \qquad (5)$$

Note that $wv^T$ is one-rank matrix. By using the Shermann-Morrison formula, we have

$$\left[ \bar{J}(U) \right]^{-1} = L^{-1} - \frac{L^{-1}wv^T L^{-1}}{1 + v^T L^{-1} w}. \qquad (6)$$

In this way, we obtain the following iterative formula

$$U^{(k+1)} = \left[ \bar{J}\left(U^{(k)}\right) \right]^{-1} \left[ b + N^{(2)}\left(U^{(k)}\right) + 2N^{(3)}\left(U^{(k)}\right) \right]$$
$$= \left( L^{-1} - \frac{L^{-1}w^{(k)}\left(v^{(k)}\right)^T L^{-1}}{1 + \left(v^{(k)}\right)^T L^{-1} w^{(k)}} \right) \left[ b + N^{(2)}\left(U^{(k)}\right) + 2N^{(3)}\left(U^{(k)}\right) \right] \qquad (7)$$

from equation (5). It is noted that only one inverse operation is required in the first step. The computing effort in one successive step is $3n^2$ multiplications, which is equivalent to that in the quasi-Newton method. To differentiate them, we name the present iterative formula as the pseudo-Newton method. If the original Jacobian matrix rather than $\bar{J}(U)$ is substituted into equation (84), we will get the original Newton method. This shows that the presented pseudo-Newton method is theoretically closely related to the original Newton method. In the above procedure, the linear coefficient matrix was used as the full-rank base matrix. We can also use the original Jacobian matrix in the first iterative step as the full-rank basis matrix. In this case, equation (2) is modified by two sides minus the matrix-vector product of the nonlinear term Jacobian matrix in the first step and vector U, and then the pseudo-Jacobian matrix in equation (3) can accordingly be redefined, and the remaining procedures can be done in the same way.

It is possible to avoid the inverse operation even in the first step. Let
$$Q = L - \Sigma l_{ii}. \tag{8}$$
By redefining the pseudo-Jacobian matrix as
$$QU + 2N^{(2)}(U) + 3N^{(3)}(U) = \left\{ \left[ QU + 2N^{(2)}(U) + 3N^{(3)}(U) \right] \left[ \frac{1}{n}\left(U^{\circ(-1)}\right)^T \right] \right\} U \tag{9}$$

We have
$$\hat{J}(U) = \Sigma l_{ii} + \left\{ \left[ QU + 2N^{(2)}(U) + 3N^{(3)}(U) \right] \left[ \frac{1}{n}\left(U^{\circ(-1)}\right)^T \right] \right\}$$
$$= \Sigma l_{ii} + sv^T \tag{10}$$

In the same procedure of deriving iterative formula (7), we get iterative formula
$$U^{(k+1)} = \left[ \hat{J}\left(U^{(k)}\right) \right]^{-1} \left[ b + N^{(2)}\left(U^{(k)}\right) + 2N^{(3)}\left(U^{(k)}\right) \right]$$
$$= \left( \Sigma l_{ii}^{-1} - \frac{\Sigma l_{ii}^{-1} s^{(k)} \left\{v^{(k)}\right\}^T \Sigma l_{ii}^{-1}}{1 + \left\{v^{(k)}\right\}^T \Sigma l_{ii}^{-1} s^{(k)}} \right) \left[ b + N^{(2)}\left(U^{(k)}\right) + 2N^{(3)}\left(U^{(k)}\right) \right] \tag{11}$$

The component-wise representation of iterative formula (11) is.
$$U_i^{(k+1)} = \frac{1}{l_{ii}} \left( z_i^{(k)} - \frac{s_i^{(k)} \sum_{j=1}^m \frac{z_j^{(k)}}{U_j^{(k)} l_{jj}}}{m + \sum_{j=1}^m \frac{s_j^{(k)}}{U_j^{(k)} l_{jj}}} \right), \tag{12}$$

where
$$z_i^{(k)} = \left\{ b + N^{(2)}\left(U^{(k)}\right) + 2N^{(3)}\left(U^{(k)}\right) \right\}_i.$$

No inverse operation is required in the pseudo-Newton iterative formula (11). The component-wise representation (12) would be preferable for the actual computations.

It is noted that the previous pseudo-Newton iterative formulas for the nonlinear polynomial-only problems avoid not only the inverse operation but also the evaluation of the Jacobian matrix. In fact, in one dimension, the pseudo-Newton method degenerates into the original Newton method, while, in contrast, the quasi-Newton method is equivalent to the secant method for the scalar nonlinear equation.

In the following, we establish the relationship between the pseudo-Jacobian matrix and original Jacobian matrix. The following theorem was establised in [2, 3].

**Theorem 2.1**: If $N^{(m)}(U)$ and $J^{(m)}(U)$ are defined as nonlinear numerical analogue of the m order nonlinear differential operator and its corresponding Jacobian matrix, respectively, then $N^{(m)}(U) = \frac{1}{m} J^{(m)}(U) U$ is always satisfied irrespective if which numerical technique is employed to discretize.

Consider nonlinear term in equation (1), by using theorem 2.1, we have
$$J_N U = 2N^{(2)}(U) + 3N^{(3)}(U), \tag{13}$$
where $J_N$ represents the Jacobian matrix of nonlinear terms. It is observed that equation (13) can not determine the Jacobian matrix $J_N$ uniquely in more than one dimension. By multiplying $\frac{1}{n}\left(U^{\circ(-1)}\right)^T$, we get
$$J_N U \left\{\frac{1}{n}\left(U^{\circ(-1)}\right)^T\right\} = \left[2N^{(2)}(U) + 3N^{(3)}(U)\right]\left\{\frac{1}{n}\left(U^{\circ(-1)}\right)^T\right\} \tag{14}$$
$$= \hat{J}_N$$
where $\hat{J}_N$ is the pseudo-Jacobian matrix. Therefore, it is clear
$$\hat{J}_N = J_N \frac{1}{n} \begin{bmatrix} 1 & \frac{U_1}{U_2} & \cdots & \frac{U_1}{U_n} \\ \frac{U_2}{U_1} & 1 & \cdots & \frac{U_2}{U_n} \\ \vdots & \vdots & \ddots & \vdots \\ \frac{U_n}{U_1} & \frac{U_n}{U_2} & \cdots & 1 \end{bmatrix} = J_N p(U), \tag{15}$$
where p(U) is defined as the deviation matrix between the original Jacobian and pseudo-Jacobian matrices for polynomial-only problems. Similar relationship for general nonlinear system is not available. It is noted that the difference between the original Jacobian and pseudo-Jacobian matrices has not effect on the stability analysis of nonlinear initial value problems given in [4]. In the pseudo-Newton method, the pseudo-Jacobian matrix is used only as the one-rank modification of the Jacobian matrix as shown in equations (4) and (8), where the linear term is still the original Jacobian matrix. Moreover, as above-mentioned, the original Jacobian matrix in the first iterative step can be considered as the full-rank basis matrix instead of L in equations (4) and (8). That means that the coefficient matrix in the pseudo-Newton method is a sum of the original Jacobian matrix in the first step and the pseudo-Jacobian matrix of the modification Jacobian matrix in the later successive step.